\author{\Large{Damanvir Singh Binner \footnote{The author acknowledges the support of IISER Mohali for providing research facilities and fellowship.}
}}
\documentclass[12pt]{article}
\usepackage[margin=1in]{geometry}
\usepackage{amsfonts,amsmath,amssymb, amsthm}
\usepackage{relsize}

\usepackage{amscd}
\usepackage{graphicx}
\usepackage{hyperref}
\usepackage{thmtools}
\usepackage{thm-restate}
\usepackage{tikz-cd} 
\usepackage{float}
\usepackage{footnote}

\DeclareMathOperator{\NR}{NR}
\begin{document}
\theoremstyle{plain} 
\newtheorem{theorem}{Theorem}
 \newtheorem{corollary}[theorem]{Corollary} 
 \newtheorem{lemma}[theorem]{Lemma} 
 \newtheorem{proposition}[theorem]{Proposition}
 
\theoremstyle{definition}
\newtheorem{definition}[theorem]{Definition}
\newtheorem{example}[theorem]{Example}
\newtheorem{conjecture}[theorem]{Conjecture}
\theoremstyle{remark}
\newtheorem{remark}[theorem]{Remark}

\title{\Large{A Complete List of all Numbers not of the Form $ax+by$}}
\date{}
\maketitle
\begin{center}
\vspace*{-8mm}
\large{Department of Mathematics \\
Indian Institute of Science Education and Research (IISER) \\
Mohali, Punjab, India \\
damanvirsinghbinner@gmail.com}
\end{center}
\begin{abstract}
For given coprime positive integers $a$ and $b$, the classical Frobenius coin problem asked to find the largest number that cannot be expressed in the form $ax+by$ for nonnegative integers $x$ and $y$, also known as the Frobenius number. Sylvester answered this question and also discovered the number of these nonrepresentable numbers. In recent times, a lot of progress has been made regarding the sums of a fixed power of these nonrepresentable numbers, commonly known as the Sylvester sums. However, the actual list of nonrepresentable numbers has remained mysterious so far. In this note, we obtain a complete list, that is a simple explicit description of these nonrepresentable numbers. We give two different proofs of our result. One of these is a short direct proof whereas the other one enables us to use Eisenstein's lattice point counting techniques to study the set of nonrepresentable numbers.
\end{abstract}
 
\section {Introduction}
\label{Intro}
Let $a$ and $b$ be given coprime positive integers $a$ and $b$. It is well known that there are finitely many natural numbers that cannot be expressed in the form $ax+by$ for nonnegative integers $x$ and $y$. These numbers are sometimes known as the \emph{nonrepresentable numbers}. The largest such number is known as the \emph{Frobenius number} of $a$ and $b$. In $1884$, Sylvester proved that the Frobenius number is $ab-a-b$, as shown in the following result, which we refer to as \emph{Sylvester's lemma}. 

\begin{lemma}[Sylvester (1882)]
\label{Frobenius}
For natural numbers $a$ and $b$ such that $\gcd(a,b)=1$, the equation
$ax+by=n$ has a solution $(x, y)$, with $x$ and $y$ nonnegative
integers, whenever $n \geq (a-1)(b-1)$. 
\end{lemma} 

Sylvester further proved there is no solution of $ax+by=n$ for exactly half of the numbers till the Frobenius number. We state this in the following theorem, which we refer to as \emph{Sylvester's theorem}.

 \begin{theorem}[Sylvester (1882)]
 \label{Sylvester}
 If $a$ and $b$ are coprime numbers, the number of natural numbers 
 that cannot be expressed in the form $ax + by$  for nonnegative integers $x$ and $y$ is equal to $\frac{(a-1)(b-1)}{2}$. 
\end{theorem}

Several curious properties of nonrepresentable numbers have been discovered in the recent years. For example, Brown and Shiue \cite{BS} discovered a formula for the sum of nonrepresentable numbers. Rødseth \cite{RS} considered the sums of higher powers of nonrepresentable numbers, commonly known as the \emph{Sylvester sums}.  A lot of progress has been made on Sylvester sums. For details, refer to \cite{B21, BThesis, komatsu2, Tu, Wang}. 

However, despite the progress on Sylvester sums, the actual list of nonrepresentable numbers seems deeply mysterious, and is the main goal of this note. 

Recently, the present author \cite[Section 3]{Binner} proved that a special case of Sylvester's theorem is equivalent to the following celebrated result of Gauss related to the law of quadratic reciprocity. 

\begin{theorem}[Gauss (1808)]
\label{Quad}
For distinct odd primes $p$ and $q$, $$\sum_{i=1}^{\frac{p-1}{2}}\Big\lfloor \frac{iq}{p} \Big\rfloor + \sum_{i=1}^{\frac{q-1}{2}}\Big\lfloor \frac{ip}{q}\Big\rfloor  =  \frac{(p-1)(q-1)}{4}, $$
\end{theorem}

Eisenstein (see \cite{Eisenstein} or \cite[Page 20]{Classical}) gave a beautiful geometric proof of Theorem $\ref{Quad}$ by using lattice point counting methods.

Further, the present author \cite[Section 2]{B21} showed that a mild generalization of Theorem \ref{Quad} is equivalent to Sylvester's theorem. One implication of the equivalence is that using Eisenstein's lattice point counting argument (used in the geometric proof of Theorem \ref{Quad}), one can obtain the number of nonrepresentable numbers. However, since the present author's proof of the equivalence as well as Eisenstein's lattice point arguments involve only counting techniques, it seems that this procedure might also yield the complete set of nonrepresentable numbers, which is usually denoted by $\NR(a,b)$. This is the goal of this note. That is, we obtain the complete list of nonrepresentable numbers, as described in Theorem \ref{Main} below. Let $R(a,b)$ denotes the set of numbers less than or equal to $a \left \lfloor\frac{b}{2} \right \rfloor+ b \left \lfloor\frac{a}{2} \right \rfloor$ that can be expressed in the form $ax+by$ for nonnegative integers $x$ and $y$.

\begin{theorem}
\label{Main}
The set $R(a,b)$ is given as 
\begin{align*}
 R(a,b) =& \left\{ai+bj: 0 \leq i \leq \Big \lfloor\frac{b}{2} \Big \rfloor, 0 \leq j \leq \Big \lfloor\frac{a}{2} \Big \rfloor \right\} \\
 & \cup \left\{a \Big \lfloor\frac{b}{2} \Big \rfloor + b \Big \lfloor\frac{a}{2} \Big \rfloor - |ai-bj|: 1 \leq i \leq \Big \lfloor\frac{b}{2} \Big \rfloor, 1 \leq j \leq \Big \lfloor\frac{a}{2} \Big \rfloor \right\},
 \end{align*}
and
 \begin{equation*}
 \NR(a,b) = \left\{0,1,2, \cdots, a \Big \lfloor\frac{b}{2} \Big \rfloor+ b \Big \lfloor\frac{a}{2} \Big \rfloor \right\} \setminus R(a,b).
 \end{equation*} 
\end{theorem}

\begin{remark}
Sylvester's theorem immediately follows from Theorem \ref{Main}. To see this, note that $$|R(a,b)| = \left( \Big \lfloor \frac{a}{2}  \Big \rfloor + 1 \right) \left( \Big \lfloor \frac{b}{2}  \Big \rfloor + 1 \right) + \Big \lfloor \frac{a}{2}  \Big \rfloor  \Big \lfloor \frac{b}{2}  \Big \rfloor.$$ Then, $$ |\NR(a,b)| =    \left( a \Big \lfloor\frac{b}{2} \Big \rfloor+ b \Big \lfloor\frac{a}{2} \Big \rfloor + 1 \right) - |R(a,b)|. $$ From here, it is straightforward to verify that $|\NR(a,b)| = \frac{(a-1)(b-1)}{2}$. 
\end{remark}

This very simple description of nonrepresentable numbers suggests that one should be able to give a short direct proof. We give such a proof in Section \ref{Short}. The readers interested in a direct proof might skip Section \ref{Proof} and directly go to Section \ref{Short}.   
Though many steps in our proof of Theorem \ref{Main} in Section \ref{Proof} are very similar to those in the proof of the equivalence in \cite[Section 3]{Binner}, we provide all the details here for the sake of completeness.
 
 \section{Proof of Theorem \ref{Main}}
 \label{Proof}
 In this section, we prove Theorem \ref{Main}. Our approach involves finding nonnegative integer solutions of $$ax+by+z =   a \Big \lfloor\frac{b}{2} \Big \rfloor+ b \Big \lfloor\frac{a}{2} \Big \rfloor. $$ Let $T$ denote the set of nonnegative integer solutions of the above equation. Then, it is easy to verify that $R(a,b)$ is the same set as
 \begin{equation}
 \label{Refer}
 \left\{ax+by: \left(x,y,a \Big \lfloor\frac{b}{2} \Big \rfloor+ b \Big \lfloor\frac{a}{2} \Big \rfloor-ax-by \right) \in T \right\}. 
 \end{equation}
 
  Thus, we basically need to study the set $T$. We recall some notation, maps and approach from \cite[Section 3]{Binner}.
 
 Let $X$, $Y$, and $Z$ denote $  \left \lfloor \frac{b}{2} \right \rfloor -x$, $  \left \lfloor \frac{a}{2}  \right \rfloor -y$, and $ b \left \lfloor \frac{a}{2} \right \rfloor -z$, respectively. Then the given equation can be rewritten as $aX + bY + Z =  b \left \lfloor \frac{a}{2} \right \rfloor$. We split our calculation into four different cases according to
\begin{enumerate}
\item $ X \geq 0, Y \geq 0, Z \geq 0$,
\item $ X \geq 0, Y \geq 0, Z < 0$,
\item $ X \geq 0, Y < 0$, or
\item $ X < 0$.
\end{enumerate}

Recall the following sets:

\begin{itemize}
\item $S_1$, $S_2$, $S_3$, and $S_4$ denote the set of nonnegative integer solutions of $ax+by+z =   a \left \lfloor\frac{b}{2} \right \rfloor+ b \left \lfloor\frac{a}{2} \right \rfloor$ in Cases $1$, $2$, $3$, and $4$, respectively.
\item $T_1$ denotes the set of nonnegative integer solutions of $ax + by + z = b \left \lfloor\frac{a}{2} \right \rfloor$.
\item $T_2$ denotes the set of nonnegative integer solutions of $ax + by + z = a \left \lfloor\frac{b}{2} \right \rfloor$. 
\item $U$ denotes the set of solutions in $T_2$ that satisfy $z=0$.
\item $V$ denotes the set of solutions in $T_2$ that satisfy $y=0$.
\item $W$ denotes the set of solutions in $T_1$ that satisfy $x=0$.
\end{itemize}

Next, recall the following maps from the sets $S_1$, $S_2$, $S_3$, and $S_4$  to $T_1$ and $T_2$. 

\begin{itemize}
\item Define $\phi_1: S_1 \rightarrow T_1$ such that $(x, y, z)\mapsto (X,Y,Z)$.
\item Define $\phi_2: S_2 \rightarrow T_2$ such that $ (x, y, z)\mapsto (x, y,-Z)$.
\item Define $\phi_3: S_3 \rightarrow T_2$ such that $(x, y, z)\mapsto (x,-Y, z)$.
\item Define $\phi_4: S_4 \rightarrow T_1$ such that $ (x, y, z) \mapsto (-X, y, z)$.
\end{itemize}

As mentioned in \cite[Section 3]{Binner}, it is easy to verify that $\phi_1$, $\phi_2$, $\phi_3$, and $\phi_4$ are well-defined injective maps and their images are given as follows:

\begin{itemize}
\item $\phi_1(S_1) = T_1$.
\item $\phi_2(S_2) = T_2 \setminus U$.
\item $\phi_3(S_3) = T_2 \setminus V$.
\item $\phi_4(S_4) = T_1 \setminus W$.
\end{itemize}

In \cite[Section 3]{Binner}, we were concerned with the cardinalities of these image sets. Presently, we note that these image sets are also easy to write out explicitly. We have 

\begin{itemize}
\item $T_1 = \left\{\left(x,y,b \left \lfloor\frac{a}{2} \right \rfloor - ax-by \right): 0 \leq y \leq \left \lfloor \frac{a}{2} \right \rfloor, 0 \leq x \leq \left \lfloor \frac{b}{a} \left(\left \lfloor \frac{a}{2} \right \rfloor - y \right) \right \rfloor \right\}$.
\item $T_2 = \left\{\left(x,y,a \left \lfloor\frac{b}{2} \right \rfloor - ax-by  \right): 0 \leq x \leq \left \lfloor \frac{b}{2} \right \rfloor, 0 \leq y \leq \left \lfloor \frac{a}{b} \left(\left \lfloor \frac{b}{2} \right \rfloor - x \right) \right \rfloor \right\}$.
\item $U = \left(\left \lfloor\frac{b}{2} \right \rfloor, 0, 0 \right)$.
\item $V = \left\{\left(x,0,a \left(\left \lfloor\frac{b}{2} \right \rfloor - x \right)\right): 0 \leq x \leq \left \lfloor \frac{b}{2} \right \rfloor \right\}$.
\item $W = \left\{\left(0,y,b \left(\left \lfloor\frac{a}{2} \right \rfloor - y \right)\right): 0 \leq y \leq \left \lfloor \frac{a}{2} \right \rfloor \right\}$.
\end{itemize}

Thus we can write the image sets as follows.

\begin{itemize}
\item $T_1 = \left\{\left(x,y,b \left \lfloor\frac{a}{2} \right \rfloor - ax-by \right): 0 \leq y \leq \left \lfloor \frac{a}{2} \right \rfloor, 0 \leq x \leq \left \lfloor \frac{b}{a} \left(\left \lfloor \frac{a}{2} \right \rfloor - y \right) \right \rfloor \right\}$.
\item $T_2 \setminus U = \left\{\left(x,y,a \left \lfloor\frac{b}{2} \right \rfloor - ax-by  \right): 0 \leq x < \left \lfloor \frac{b}{2} \right \rfloor, 0 \leq y \leq \left \lfloor \frac{a}{b} \left(\left \lfloor \frac{b}{2} \right \rfloor - x \right) \right \rfloor \right\}$.
\item $T_2 \setminus V = \left\{\left(x,y,a \left \lfloor\frac{b}{2} \right \rfloor - ax-by  \right): 0 \leq x \leq \left \lfloor \frac{b}{2} \right \rfloor, 0 < y \leq \left \lfloor \frac{a}{b} \left(\left \lfloor \frac{b}{2} \right \rfloor - x \right) \right \rfloor \right\}$.
\item $T_1 \setminus W = \left\{\left(x,y,b \left \lfloor\frac{a}{2} \right \rfloor - ax-by \right): 0 \leq y \leq \left \lfloor \frac{a}{2} \right \rfloor, 0 < x \leq \left \lfloor \frac{b}{a} \left(\left \lfloor \frac{a}{2} \right \rfloor - y \right) \right \rfloor \right\}$.
\end{itemize}

Since these sets are in bijection with the sets $S_1$, $S_2$, $S_3$ and $S_4$ via the maps $\phi_1$, $\phi_2$, $\phi_3$, and $\phi_4$, which are very simple bijections, we get that the sets $S_1$, $S_2$, $S_3$ and $S_4$ are given as follows.

\begin{itemize}
\item $S_1 = \left\{\left( \left \lfloor\frac{b}{2} \right \rfloor - i, j, b\left \lfloor\frac{a}{2} \right \rfloor + ai-bj \right): 0 \leq j \leq \left \lfloor\frac{a}{2} \right \rfloor, 0 \leq i \leq \left \lfloor \frac{bj}{a} \right \rfloor \right\}$.
\item $S_2 = \left\{\left( \left \lfloor\frac{b}{2} \right \rfloor + i, \left \lfloor\frac{a}{2} \right \rfloor - j, bj-ai \right): 0 \leq j \leq \left \lfloor\frac{a}{2} \right \rfloor, 0 < i \leq \left \lfloor \frac{bj}{a} \right \rfloor \right\}$.
\item $S_3 = \left\{\left( \left \lfloor\frac{b}{2} \right \rfloor - i, j, b\left \lfloor\frac{a}{2} \right \rfloor + ai-bj \right): 0 < i \leq \left \lfloor\frac{b}{2} \right \rfloor, 0 \leq j \leq \left \lfloor \frac{ai}{b} \right \rfloor \right\}$.
\item $S_4 = \left\{\left( \left \lfloor\frac{b}{2} \right \rfloor - i, \left \lfloor\frac{a}{2} \right \rfloor + j, ai-bj \right): 0 \leq i \leq \left \lfloor\frac{b}{2} \right \rfloor, 0 < j \leq \left \lfloor \frac{ai}{b} \right \rfloor \right\}$.
\end{itemize}

The union of the sets $S_i$ then gives the set $T$. Therefore, by \eqref{Refer}, we get that the set $R(a,b)$ is given as

\begin{equation}
 \label{Final}
  A_1 \cup A_2 \cup A_3 \cup A_4, 
 \end{equation}
 where the sets $A_i$ are as follows.
 
 \begin{itemize}
\item $A_1 = \left\{a\left( \left \lfloor\frac{b}{2} \right \rfloor - i \right) + bj: 0 \leq j \leq \left \lfloor\frac{a}{2} \right \rfloor, 0 \leq i \leq \left \lfloor \frac{bj}{a} \right \rfloor \right\}$.
\item $A_2 = \left\{a\left( \left \lfloor\frac{b}{2} \right \rfloor + i \right) + b \left(\left \lfloor\frac{a}{2} \right \rfloor - j  \right): 0 \leq j \leq \left \lfloor\frac{a}{2} \right \rfloor, 0 < i \leq \left \lfloor \frac{bj}{a} \right \rfloor \right\}$.
\item $A_3 = \left\{a\left( \left \lfloor\frac{b}{2} \right \rfloor - i \right) + bj: 0 < i \leq \left \lfloor\frac{b}{2} \right \rfloor, 0 \leq j \leq \left \lfloor \frac{ai}{b} \right \rfloor \right\}$.
\item $A_4 = \left\{a\left( \left \lfloor\frac{b}{2} \right \rfloor - i \right) + b \left( \left \lfloor\frac{a}{2} \right \rfloor + j \right): 0 \leq i \leq \left \lfloor\frac{b}{2} \right \rfloor, 0 < j \leq \left \lfloor \frac{ai}{b} \right \rfloor \right\}$.
\end{itemize}

Using Eisenstein's lattice point argument in the geometric proof of Theorem \ref{Quad}, we can combine $A_1$ and $A_3$ to get $$ A_1 \cup A_3 = \left\{ a \left \lfloor \frac{b}{2} \right \rfloor - (ai-bj): 0 \leq i \leq \left \lfloor \frac{b}{2} \right \rfloor, 0 \leq j \leq \left \lfloor \frac{a}{2} \right \rfloor \right\}. $$ Similarly, we can combine $A_2$ and $A_4$ to get $$ A_2 \cup A_4 = \left\{ a \left \lfloor \frac{b}{2} \right \rfloor + b \left \lfloor \frac{a}{2} \right \rfloor - |ai-bj|: 1 \leq i \leq \left \lfloor \frac{b}{2} \right \rfloor, 1 \leq j \leq \left \lfloor \frac{a}{2} \right \rfloor \right\}. $$ Thus, the set of numbers less than or equal to $a \left \lfloor\frac{b}{2} \right \rfloor+ b \left \lfloor\frac{a}{2} \right \rfloor$ that can be expressed in the form $ax+by$ is given by 
\begin{equation*}
\begin{aligned}
 R(a,b) =& A_1 \cup A_2 \cup A_3 \cup A_4 \\
 =& \left\{ai: 0 \leq i \leq \Big \lfloor\frac{b}{2} \Big \rfloor \right\}  \\
 & \cup \left\{a \Big \lfloor\frac{b}{2} \Big \rfloor + bj: 1 \leq j \leq \Big \lfloor\frac{a}{2} \Big \rfloor \right\}  \\
 & \cup \left\{a \Big \lfloor\frac{b}{2} \Big \rfloor - (ai-bj): 1 \leq i \leq \Big \lfloor\frac{b}{2} \Big \rfloor, 1 \leq j \leq \Big \lfloor\frac{a}{2} \Big \rfloor \right\}  \\
 & \cup \left\{a \Big \lfloor\frac{b}{2} \Big \rfloor + b \Big \lfloor\frac{a}{2} \Big \rfloor - |ai-bj|: 1 \leq i \leq \Big \lfloor\frac{b}{2} \Big \rfloor, 1 \leq j \leq \Big \lfloor\frac{a}{2} \Big \rfloor \right\}.
 \end{aligned}
\end{equation*}
We can rewrite the third set in the above union as: $$ \left\{ai+bj: 0 \leq i < \Big \lfloor\frac{b}{2} \Big \rfloor, 1 \leq j \leq \Big \lfloor\frac{a}{2} \Big \rfloor \right\} .$$ Then, we can combine the first three sets to get the set $$ \left\{ai+bj: 0 \leq i \leq \Big \lfloor\frac{b}{2} \Big \rfloor, 0 \leq j \leq \Big \lfloor\frac{a}{2} \Big \rfloor \right\} .$$ Therefore, we can write $R(a,b)$ as 

\begin{align*}
 R(a,b) =& \left\{ai+bj: 0 \leq i \leq \Big \lfloor\frac{b}{2} \Big \rfloor, 0 \leq j \leq \Big \lfloor\frac{a}{2} \Big \rfloor \right\} \\
 & \cup \left\{a \Big \lfloor\frac{b}{2} \Big \rfloor + b \Big \lfloor\frac{a}{2} \Big \rfloor - |ai-bj|: 1 \leq i \leq \Big \lfloor\frac{b}{2} \Big \rfloor, 1 \leq j \leq \Big \lfloor\frac{a}{2} \Big \rfloor \right\}.
 \end{align*}

By Sylvester's lemma, $a \left \lfloor\frac{b}{2} \right \rfloor+ b \left \lfloor\frac{a}{2} \right \rfloor$ is a number greater than the Frobenius number of $a$ and $b$. Therefore, the set $\NR(a,b)$ of numbers that cannot be expressed in the form $ax+by$ can be obtained by removing the elements of $R(a,b)$ from the set of numbers in the interval $[0, a \left \lfloor\frac{b}{2} \right \rfloor+ b \left \lfloor\frac{a}{2} \right \rfloor]$. That is, we have 
\begin{equation*}
 \NR(a,b) = \left\{0,1,2, \cdots, a \Big \lfloor\frac{b}{2} \Big \rfloor+ b \Big \lfloor\frac{a}{2} \Big \rfloor \right\} \setminus R(a,b). 
 \end{equation*} 
 
 This completes the proof of Theorem \ref{Main}.
 
 \section{A short direct proof of Theorem \ref{Main}}
\label{Short}

We have $$R(a,b) = \left\{ax + by : x \geq 0, y \geq 0, ax+by \leq a \Big \lfloor\frac{b}{2} \Big \rfloor+ b \Big \lfloor\frac{a}{2} \Big \rfloor \right\}. $$ Note that for $ax+by \leq a \Big \lfloor\frac{b}{2} \Big \rfloor+ b \Big \lfloor\frac{a}{2} \Big \rfloor$, it is not possible to have $x > \left \lfloor\frac{b}{2} \right \rfloor$ and $y > \left \lfloor\frac{a}{2} \right \rfloor$. Thus, we have the following three cases.

\begin{enumerate}
\item $x \leq \left \lfloor\frac{b}{2} \right \rfloor$, $y \leq \left \lfloor\frac{a}{2} \right \rfloor$.
\item $x >  \left \lfloor\frac{b}{2} \right \rfloor$, $y < \left \lfloor\frac{a}{2} \right \rfloor$.
\item $y > \left \lfloor\frac{a}{2} \right \rfloor$, $x < \left \lfloor\frac{b}{2} \right \rfloor$.
\end{enumerate} 

Next, we find the contribution to $R(a,b)$ from each of these three cases. 

Case $1$: Let $x \leq \left \lfloor\frac{b}{2} \right \rfloor$ and $y \leq \left \lfloor\frac{a}{2} \right \rfloor$. In this case, it is clear that $ax+by \leq a \left \lfloor\frac{b}{2} \right \rfloor+ b \left \lfloor\frac{a}{2} \right \rfloor$. Thus the contribution to $R(a,b)$ from Case $1$ is given by $$ \left\{ax+by: 0 \leq x \leq \Big \lfloor\frac{b}{2} \Big \rfloor, 0 \leq y \leq \Big \lfloor\frac{a}{2} \Big \rfloor \right\}. $$

Case $2$: Let $x >  \left \lfloor\frac{b}{2} \right \rfloor$ and $y < \left \lfloor\frac{a}{2} \right \rfloor$. Set $X = x-\left \lfloor\frac{b}{2} \right \rfloor $ and $Y = \left \lfloor\frac{a}{2} \right \rfloor - y $. Note that $X \geq 1$ and $ 1 \leq Y \leq \left \lfloor\frac{a}{2} \right \rfloor$. Then the condition $ax+by \leq a \left \lfloor\frac{b}{2} \right \rfloor+ b \left \lfloor\frac{a}{2} \right \rfloor$ can be rewritten as $aX \leq bY$. Therefore, $$ 1 \leq X \leq \left \lfloor \frac{bY}{a} \right \rfloor \leq \left \lfloor \frac{b \left \lfloor \frac{a}{2}\right \rfloor}{a} \right \rfloor \leq \left \lfloor \frac{b}{2}\right \rfloor. $$ Then, 
\begin{align*}
 ax + by &= a\left(X + \left \lfloor\frac{b}{2} \right \rfloor \right) + b\left(\left \lfloor\frac{a}{2} \right \rfloor - Y \right) \\
 &= a \Big \lfloor\frac{b}{2} \Big \rfloor+ b \Big \lfloor\frac{a}{2} \Big \rfloor + (aX-bY). 
 \end{align*}
 
 Thus, the contribution to $R(a,b)$ from Case $2$ is given by $$ \left\{ a \Big \lfloor\frac{b}{2} \Big \rfloor+ b \Big \lfloor\frac{a}{2} \Big \rfloor + (aX-bY): 1 \leq X \leq \Big \lfloor\frac{b}{2} \Big \rfloor, 1 \leq Y \leq \Big \lfloor\frac{a}{2} \Big \rfloor, aX \leq bY  \right\}. $$ Note that it is clear that the members of the above set indeed lie in $R(a,b)$. 
 A similar analysis shows that the contribution to $R(a,b)$ from Case $3$ is given by $$ \left\{ a \Big \lfloor\frac{b}{2} \Big \rfloor+ b \Big \lfloor\frac{a}{2} \Big \rfloor - (aX-bY): 1 \leq X \leq \Big \lfloor\frac{b}{2} \Big \rfloor, 1 \leq Y \leq \Big \lfloor\frac{a}{2} \Big \rfloor, aX \geq bY  \right\}. $$ Combining the contributions from Cases $2$ and $3$, we get that the net contribution to $R(a,b)$ from these cases is given by $$ \left\{a \Big \lfloor\frac{b}{2} \Big \rfloor + b \Big \lfloor\frac{a}{2} \Big \rfloor - |aX-bY|: 1 \leq X \leq \Big \lfloor\frac{b}{2} \Big \rfloor, 1 \leq Y \leq \Big \lfloor\frac{a}{2} \Big \rfloor \right\}. $$ Combining this with the contribution to $R(a,b)$ from Case $1$ completes the proof of first part of Theorem \ref{Main}. The second part of Theorem \ref{Main} simply follows from the fact that $a \left \lfloor\frac{b}{2} \right \rfloor+ b \left \lfloor\frac{a}{2} \right \rfloor$ is a number greater than the Frobenius number of $a$ and $b$, which is a consequence of Sylvester's lemma.
 
 \section{Examples}
\label{Calc}
Let $a=7$ and $b=5$. By Sylvester's lemma and Sylvester's theorem, we know that every number greater than or equal to $24$ can be expressed in the form $7x+5y$ and exactly $12$ numbers are nonrepresentable. Here, we find all the nonrepresentable numbers using Theorem \ref{Main}. We have $$ \NR(7,5) = \left\{0,1,2, \cdots, 29 \right\} \setminus R(7,5), $$
where
\begin{equation}
\label{List24}
\begin{aligned}
 R(7,5) =& \left\{7i+5j: 0 \leq i \leq 2, 0 \leq j \leq 3 \right\} \\
 & \cup \left\{ 29 - |7i-5j|: 1 \leq i \leq 2, 1 \leq j \leq 3 \right\}.
 \end{aligned}
 \end{equation}

In the first set in \eqref{List24}, the values when either $i$ or $j$ is $0$ are given by $\{0\} \cup \left\{7i: 1 \leq i \leq 2 \right\} \cup  \left\{5j: 1 \leq j \leq 3 \right\} = \{0\} \cup \{7,14\} \cup \{5,10,15\} = \{0,5,7,10,14,15\}. $   For finding the remaining values in the two sets in \eqref{List24}, we construct Table \ref{75}.

\begin{table}[htpb]
    \centering
    \begin{tabular}{|c|c|c|c|c|c|}
 \hline
$i$  & $j$  & $7i+5j$ & $7i-5j$ & $|7i-5j|$ & $29 - |7i-5j|$ \\
\hline
1 & 1 & 12 & 2 & 2 & 27 \\
\hline 
1 & 2 & 17 & -3 & 3 & 26 \\
\hline 
1 & 3 & 22 & -8 & 8 & 21 \\
\hline
2 & 1 & 19 & 9 & 9 & 20 \\
\hline 
2 & 2 & 24 & 4 & 4 & 25 \\
\hline 
2 & 3 & 29 & -1 & 1 & 28 \\
\hline
\end{tabular}
\caption{The calculation of $\NR(7,5)$.}
        \label{75}
\end{table}

Thus, we have 
\begin{align*}
 R(7,5) &=  \{0,5,7,10,14,15\} \\
  & \cup \left\{12,17,22,19,24,29 \right\} \\
 & \cup \left\{ 27,26,21,20,25,28 \right\}.
 \end{align*}
 
 Therefore, $$ R(7,5) = \left\{ 0,5,7,10,12,14,15,17,19,20,21,22,24,25,26,27,28,29 \right\}, $$ and thus 
 \begin{align*}
 \NR(7,5) &= \left\{0,1,2, \cdots, 29 \right\} \setminus R(7,5) \\
 &= \left\{0,1,2, \cdots, 29 \right\} \setminus \left\{ 0,5,7,10,12,14,15,17,19,20,21,22,24,25,26,27,28,29 \right\} \\
 &= \left\{ 1,2,3,4,6,8,9,11,13,16,18,23 \right\}. 
 \end{align*}

We describe our method for another much lengthier example. Let $a=29$ and $b=23$. By Sylvester's lemma and Sylvester's theorem, we know that every number greater than or equal to $616$ can be expressed in the form $29x+23y$ and exactly $308$ numbers are nonrepresentable. 
Here, we find all the nonrepresentable numbers using Theorem \ref{Main}. We have $$ \NR(29,23) = \left\{0,1,2, \cdots, 641 \right\} \setminus R(29,23), $$
where
\begin{equation}
\label{List}
\begin{aligned}
 R(29,23) &= \left\{29i+23j: 0 \leq i \leq 11, 0 \leq j \leq 14 \right\} \\
 & \cup \left\{ 641 - |29i-23j|: 1 \leq i \leq 11, 1 \leq j \leq 14 \right\}.
 \end{aligned}
 \end{equation}

In the first set in \eqref{List}, the values when either $i$ or $j$ is $0$ are given by $\{0\} \cup \left\{29i: 1 \leq i \leq 11 \right\} \cup  \left\{23j: 1 \leq j \leq 14 \right\}$, which can be written explicitly as 
\begin{equation}
\label{List1}
\begin{aligned}
 \{0\} & \cup \{29,59,87,116,145,174,203,232,261,290,319\} \\
 & \cup \{23,46,69,92,115,138,161,184,207,230,253,276,299,322\}.
\end{aligned}
 \end{equation}
  
The calculation for finding the remaining values in the two sets in \eqref{List} is somewhat lengthy and it is time saving to work for a fixed value of $j$ and all values of $i$. As an example, we demonstrate the computation of these sets in the case $j=4$ and $ 1 \leq i \leq 11$. Other cases are similar. For $j=4$, we have $29i + 23j = 92+29i$, and its set of values as $1 \leq i \leq 11$ is given as $$ \{121, 150, 179, 208, 237, 266, 295, 324, 353, 382, 411 \}. $$ Similarly, for $j=4$, we have $29i - 23j = 29i-92$, and its set of values as $1 \leq i \leq 11$ is given as $$\{-63,-34,-5,24,53,82,111,140, 169,198,227\},$$ and thus the set of values of $|29i-23j|$ is given as $$\{63,34,5,24,53,82,111,140, 169,198,227\}.$$ Therefore, the set of values of $641-|29i-23j|$ is given as $$\{578, 607, 636, 617,588, 559, 530, 501, 472, 443, 414\}.$$ 

The calculation for other values of $j$ can be done with similar ease, and the answers are recorded in Table \ref{2923}.

\begin{table}[htpb]
    \centering
    \begin{tabular}{|c|c|c|}
 \hline
$j$ & $29i+23j: 1 \leq i \leq 11$ & $641-|29i-23j|: 1 \leq i \leq 11$  \\
\hline
1 & $52, 81, 110, 139, 168, 197, 226, 255, 284, 313, 342$  & $635, 606, 577, 548, 519, 490, 461, 432, 403, 374, 345$  \\
\hline 
2 & $75, 104, 133, 162, 191, 220, 249, 278, 307, 336, 365$ & $624, 629, 600, 571, 542, 513, 484, 455, 426, 397, 368$ \\
\hline 
3 & $98, 127, 156, 185, 214, 243, 272, 301, 330, 359, 388$  & $601, 630, 623, 594, 565, 536, 507, 478, 449, 420, 391$ \\
\hline
4 & $121, 150, 179, 208, 237, 266, 295, 324, 353, 382, 411 $ & $578, 607, 636, 617,588, 559, 530, 501, 472, 443, 414$ \\
\hline 
5 & $144, 173, 202, 231, 260, 289, 318, 347, 376, 405, 434$ & $555, 584, 613, 640, 611, 582, 553, 524, 495, 466, 437$ \\
\hline 
6 & $167, 196, 225, 254, 283, 312, 341, 370, 399, 428, 457$ & $532, 561, 590, 619, 634, 605, 576, 547, 518, 489, 460$ \\
\hline 
7 & $190, 219, 248, 277, 306, 335, 364, 393, 422, 451, 480$ & $509, 538, 567, 596, 625, 628, 599, 570, 541, 512, 483$ \\
\hline 
8 & $213, 242, 271, 300, 329, 358, 387, 416, 445, 474, 503$ & $486, 515, 544, 573, 602, 631, 622, 593, 564, 535, 506$  \\
\hline 
9 & $236, 265, 294, 323, 352, 381, 401, 439, 468, 497, 526$ & $463, 492, 521, 550, 579, 608, 637, 616, 587, 558, 529$ \\
\hline 
10 & $259, 288, 317, 346, 375, 404, 433, 462, 491, 520, 549$ &  $440, 469, 498, 527, 556, 585, 614, 639, 610, 581, 552$ \\
\hline 
11 & $282, 311, 340, 369, 398, 427, 456, 485, 514, 543, 572$ & $417, 446, 475, 504, 533, 562, 591, 620, 633, 604, 575$  \\
\hline 
12 & $305, 334, 363, 392, 421, 450, 479, 508, 537, 566, 595$ & $394, 423, 452, 481, 510, 539, 568, 597, 626, 627, 598$ \\
\hline 
13 & $328, 357, 386, 415, 444, 473, 502, 531, 560, 589, 618$ & $371, 400, 429, 458, 487, 516, 545, 574, 603, 632, 621$  \\
\hline 
14 & $351, 380, 409, 438, 467, 496, 525, 554, 583, 612, 641$ & $348, 377, 406, 435, 464, 493, 522, 551, 580, 609, 638$ \\
\hline
\end{tabular}
\caption{The calculation of $\NR(29,23)$.}
        \label{2923}
\end{table}

Collecting all the values in Table \ref{2923}, along with those in \eqref{List1} and then sorting, we get that $R(29,23)$ consists of the following numbers.

\begin{equation}
\label{List2}
\begin{aligned}
& 0, 23, 29, 46, 52, 58, 69, 75, 81, 87, 92, 98, 104, 110, 115, 116, 121, 127, \\
& 133, 138, 139, 144, 145, 150, 156, 161, 162, 167, 168, 173, 174, 179, 184, 185, \\
& 190, 191, 196, 197, 202, 203, 207, 208, 213, 214, 219, 220, 225, 226, 230, 231, \\
&232, 236, 237, 242, 243, 248, 249, 253, 254, 255, 259, 260, 261, 265, 266, 271, \\
& 272, 276, 277, 278, 282, 283, 284, 288, 289, 290, 294, 295, 299, 300, 301, 305, \\
& 306, 307, 311, 312, 313, 317, 318, 319, 319, 322, 322, 323, 324, 328, 329, 330, \\
& 334, 335, 336, 340, 341, 342, 342, 345, 346, 347, 348, 351, 351, 352, 353, 357, \\
& 358, 359, 363, 364, 365, 365, 368, 369, 370, 371, 374, 375, 376, 377, 380, 380, \\
& 381, 382, 386, 387, 388, 388, 391, 392, 393, 394, 397, 398, 399, 400, 403, 404, \\
& 405, 406, 409, 409, 410, 411, 411, 414, 415, 416, 417, 420, 421, 422, 423, 426, \\
& 427, 428, 429, 432, 433, 434, 434, 435, 437, 438, 438, 439, 440, 443, 444, 445, \\
& 446, 449, 450, 451, 452, 455, 456, 457, 457, 458, 460, 461, 462, 463, 464, 466, \\
& 467, 467, 468, 469, 472, 473, 474, 475, 478, 479, 480, 480, 481, 483, 484, 485, \\
& 486, 487, 489, 490, 491, 492, 493, 495, 496, 496, 497, 498, 501, 502, 503, 503, \\
& 504, 506, 507, 508, 509, 510, 512, 513, 514, 515, 516, 518, 519, 520, 521, 522, \\
& 524, 525, 525, 526, 526, 527, 529, 530, 531, 532, 533, 535, 536, 537, 538, 539, \\
& 541, 542, 543, 544, 545, 547, 548, 549, 549, 550, 551, 552, 553, 554, 554, 555, \\
& 556, 558, 559, 560, 561, 562, 564, 565, 566, 567, 568, 570, 571, 572, 572, 573, \\
& 574, 575, 576, 577, 578, 579, 580, 581, 582, 583, 583, 584, 585, 587, 588, 589, \\
& 590, 591, 593, 594, 595, 595, 596, 597, 598, 599, 600, 601, 602, 603, 604, 605, \\
& 606, 607, 608, 609, 610, 611, 612, 612, 613, 614, 616, 617, 618, 618, 619, 620, \\
& 621, 622, 623, 624, 625, 626, 627, 628, 629, 630, 631, 632, 633, 634, 635, 636, \\
& 637, 638, 639, 640, 641, 641.
\end{aligned}
 \end{equation}
 
 Removing the numbers in \eqref{List2} from the set $\{0,1,2, \cdots, 641\}$, we find the set $N(29,23)$. That is, the complete list of nonnegative integers that cannot be expressed in the form $29x+23y$ for nonnegative integers $x$ and $y$ is given as follows.
 
 \begin{equation*}
 \begin{aligned}
 & 1, 2, 3, 4, 5, 6, 7, 8, 9, 10, 11, 12, 13, 14, 15, 16, 17, 18, 19, 20, 21, 22, \\
& 24, 25, 26, 27, 28, 30, 31, 32, 33, 34, 35, 36, 37, 38, 39, 40, 41, 42, 43, 44, \\
& 45, 47, 48, 49, 50, 51, 53, 54, 55, 56, 57, 59, 60, 61, 62, 63, 64, 65, 66, 67, \\
& 68, 70, 71, 72, 73, 74, 76, 77, 78, 79, 80, 82, 83, 84, 85, 86, 88, 89, 90, 91, \\
& 93, 94, 95, 96, 97, 99, 100, 101, 102, 103, 105, 106, 107, 108, 109, 111, 112, \\
& 113, 114, 117, 118, 119, 120, 122, 123, 124, 125, 126, 128, 129, 130, 131, 132, \\
& 134, 135, 136, 137, 140, 141, 142, 143, 146, 147, 148, 149, 151, 152, 153, 154, \\
& 155, 157, 158, 159, 160, 163, 164, 165, 166, 169, 170, 171, 172, 175, 176, 177, \\
& 178, 180, 181, 182, 183, 186, 187, 188, 189, 192, 193, 194, 195, 198, 199, 200, \\
& 201, 204, 205, 206, 209, 210, 211, 212, 215, 216, 217, 218, 221, 222, 223, 224, \\
& 227, 228, 229, 233, 234, 235, 238, 239, 240, 241, 244, 245, 246, 247, 250, 251, \\
& 252, 256, 257, 258, 262, 263, 264, 267, 268, 269, 270, 273, 274, 275, 279, 280, \\
& 281, 285, 286, 287, 291, 292, 293, 296, 297, 298, 302, 303, 304, 308, 309, 310, \\
& 314, 315, 316, 320, 321, 325, 326, 327, 331, 332, 333, 337, 338, 339, 343, 344, \\
& 349, 350, 354, 355, 356, 360, 361, 362, 366, 367, 372, 373, 378, 379, 383, 384, \\
& 385, 389, 390, 395, 396, 401, 402, 407, 408, 412, 413, 418, 419, 424, 425, 430, \\
& 431, 436, 441, 442, 447, 448, 453, 454, 459, 465, 470, 471, 476, 477, 482, 488, \\
&494, 499, 500, 505, 511, 517, 523, 528, 534, 540, 546, 557, 563, 569, 586, 592, \\
& 615.
\end{aligned}
 \end{equation*}

\end{document}